\documentclass[reqno]{amsart}
\usepackage{amscd}
\usepackage{amssymb}

\newcommand{\bbC}{{\mathbb C}}
\newcommand{\bbH}{{\mathbb H}}

\newcommand{\bbP}{{\mathbb P}}
\newcommand{\bbR}{{\mathbb R}}
\newcommand{\bbZ}{{\mathbb Z}}
\newcommand{\cO}{{\mathcal O}}
\newcommand{\bbQ}{{\mathbb Q}}
\newcommand{\M}{{\mathcal M}}

\newcommand{\cE}{{\mathcal E}}
\newcommand{\C}{{\bbC}}
\newcommand{\R}{{\bbR}}
\newcommand{\<}{\langle}
\renewcommand{\>}{\rangle}

\DeclareMathOperator{\tr}{tr}

\DeclareMathOperator{\Hess}{Hess}

\newtheorem{lemma}{Lemma}[section]
\newtheorem{prop}[lemma]{Proposition}
\newtheorem{thm}{Theorem}

\newtheorem{cor}[thm]{Corollary}
\newtheorem{defn}{Definition}

\theoremstyle{plain}

\theoremstyle{definition}

\theoremstyle{remark}

\newtheorem{note}{Note}[section]
\newtheorem{remark}{Remark}[section]

\begin{document}
\title{Hyperbolic monopoles and holomorphic spheres}

\author{Michael K. Murray}
\address[Michael K. Murray]
{Department of Pure Mathematics\\
University of Adelaide\\
Adelaide, SA 5005 \\
Australia}
\email[Michael K. Murray]{mmurray@maths.adelaide.edu.au}

\author{Paul Norbury}
\address[Paul Norbury]
{Department of Pure Mathematics\\
University of Adelaide\\
Adelaide, SA 5005 \\
Australia}
\email[Paul Norbury]{pnorbury@maths.adelaide.edu.au}

\author{Michael A. Singer}
\address[Michael A. Singer]
{Department of Mathematics and Statistics \\
James Clerk Maxwell
Building \\
University of Edinburgh EH9 3JZ \\
U.K.}
\email[Michael A. Singer]{michael@maths.ed.ac.uk}

\thanks{MKM and PN acknowledge the support of the Australian
Research Council}

\keywords{}
\subjclass{81T13, 53C07}

\begin{abstract}
     We associate to an $SU(2)$ hyperbolic monopole a holomorphic
     sphere embedded in projective space and use this to uncover
     various features of the monopole.
\end{abstract}
\maketitle

\section{Introduction}
In this paper we exploit the geometry of hyperbolic space to study
monopoles.  We will use features of hyperbolic space that do not arise
in Euclidean space, and hence expose properties of hyperbolic
monopoles that have no analogues for Euclidean monopoles.  The space
of geodesics in $\bbH^3$ is the complex manifold
$$
Z=\bbP^1\times\bbP^1-\bar{\Delta}
$$
  where the point
$(w,z)\in\bbP^1\times\bbP^1$ represents the geodesic that runs from
$\hat{w}=-1/\bar{w}$, the antipodal point of $w$, to $z$ considered as
points on the sphere at infinity.  The antidiagonal $\bar{\Delta}$ has
been removed, although one aspect of this paper is that in some sense
we can replace the antidiagonal, making sense of $(\hat{z},z)$, which
represents a geodesic from $z$ to itself.

A monopole is a pair $(A,\Phi)$ consisting of a connection $A$ with
$L^2$ curvature $F_A$ defined on a trivial bundle $E$ over $\bbR^3$
with structure group $SU(2)$, and a Higgs field
$\Phi:\bbR^3\rightarrow{\bf su}(2)$ that solves the Bogomolny equation
\[ d_A\Phi=*F_A \]
and satisfies $\lim_{r\rightarrow\infty}||\Phi||=m$,
the mass of the monopole.  The charge of the monopole is defined to be
the topological degree of the map $\Phi_{\infty}:S^2\rightarrow S^2$.
The gauge group $\bbR^3\rightarrow SU(2)$ acts on the equations and we
identify gauge equivalent monopoles.  The metric is featured in the
Hodge star, $*$.  In this paper we will mainly consider the hyperbolic
metric and sometimes refer to the Euclidean metric.

Hyperbolic and Euclidean monopoles have been studied using three
constructions: the spectral curve, Nahm data, and the rational map
\cite{AtiMag,ABrBou,DonNah,HitCon,HitMon,MSiCom,MSiSpe}.  Each
construction is related, although they are different enough that a
particular aspect of monopoles is often more readily seen from the
perspective of one of these constructions.  The rational map and the
spectral curve use solutions of an ordinary differential equation
((\ref{eq:scat}) in the next section) defined along geodesics, known
as the {\em scattering equation}.  The rational map arises when one
restricts the scattering equation to the pencil of geodesics that
contain a common point (which may be at infinity.)  The spectral curve
is defined to be the set of geodesics along which the scattering
equation has an $L^2$ solution.  It is a compact algebraic curve
inside the variety of geodesics.  The spectral curve of a charge $k$
monopole is a degree $(k,k)$ curve in
$\bbP^1\times\bbP^1-\bar{\Delta}$, or equivalently the zero set of a
holomorphic section of the line bundle $\cO(k,k)$.

In this paper, we introduce a fourth construction for hyperbolic
monopoles---an embedded holomorphic sphere in projective space.
Such a construction has existed previously for half-integer mass
hyperbolic monopoles \cite{ABrBou}.  We define it for any real mass
using completely different techniques to \cite{ABrBou} and give
various applications, some of which are known for half-integer mass
monopoles.

A hyperbolic monopole has a well-defined limit at infinity given by a
reducible connection $A_{\infty}$ over a two-sphere.  Denote by
$F_{A_{\infty}}$ the curvature of the reducible connection.
\begin{thm}   \label{th:sphere}
     An $SU(2)$ hyperbolic monopole $(A,\Phi)$ of charge $k$ is
     determined by a degree $k$ holomorphic embedding
     $q:\bbP^1\rightarrow\bbP^k$ uniquely defined up to the action of
     $U(k+1)$ on its image with the properties:

     (i) $\Sigma=\{(w,z)\in\bbP^1\times\bbP^1\ |\ \langle
     q(\hat{w}),q(z)\rangle=0\}$ is the spectral curve of
     $(A,\Phi)$;

     (ii) $F_{A_{\infty}}=q^*\omega$, for $\omega$ the
     Kahler form on $\bbP^k$.
\end{thm}

One consequence of the theorem is the fact that an $SU(2)$ hyperbolic
monopole is determined up to gauge by its reducible connection on the
sphere at infinity.  This was proven in \cite{NorBou} by a
different method, and for the half-integer mass case in \cite{ABrBou}
also using an embedded sphere in projective space.

Theorem~\ref{th:sphere} relies on the fact that
$\bbP^1\times\bbP^1-\bar{\Delta}$, the twistor space of geodesics in
$\bbH^3$, has a compactification obtained by including a totally real
surface.  The same situation seems to arise for spherical monopoles
\cite{PauMon} which would lead one to predict that a monopole with one
singularity on $S^3$ is determined by its asymptotic value near the
singularity, and it is neatly described by a holomorphic sphere in
projective space.  The construction does not apply to Euclidean
monopoles.  Only the charge of a Euclidean monopole is detected from
its reducible connection at infinity.  The difference comes down to
the asymptotic decay conditions forced on finite energy monopoles in
Euclidean and hyperbolic spaces.

The {\em centre} of a Euclidean monopole is defined in \cite{AHiGeo}.
Previously a definition of the centre has existed only for
half-integer mass hyperbolic monopoles.  In \cite{SinCen} the third
author proposed a definition for a general hyperbolic monopole but
could not prove that the centre is unique.  Intuitively, the centre of
a monopole arises from the $PSL(2,\bbC)$ action on hyperbolic space.
One would like to show that the $PSL(2,\bbC)$ orbit of a hyperbolic
monopole possesses a centred monopole unique up to the action of
$SO(3)\subset PSL(2,\bbC)$.  The holomorphic sphere allows one to
apply geometric invariant theory to obtain such a definition for the
centre of the monopole.

\begin{thm}  \label{th:git}
     There is a lift of the $PSL(2,\bbC)$ action on the space of
     hyperbolic monopoles to a linear $SL(2,\bbC)$ action on $\bbC^N$
     whose stable points contain the space of hyperbolic monopoles.  A
     monopole has a unique centre, and it is centred when it lies in
     the zero set of the moment map for $SU(2)\subset SL(2,\bbC)$.
\end{thm}

Given a hyperbolic monopole $(A,\Phi)$ and a point $w\in
S^2_{\infty}$, one can use the scattering equation along geodesics
$\gamma$ satisfying $\lim_{t\rightarrow -\infty}\gamma(t)=w$ to
define a degree $k$ rational map
$f_{w}(z):\bbP^1\rightarrow\bbP^1$.  Moreover, the rational map
uniquely determines the monopole.  Previously, it has not been
understood how one might relate the different rational maps as $w$
is varied.  The holomorphic sphere $q$ in some sense combines all of
these rational maps.
\begin{thm}  \label{th:rat}
     For any $w\in S^2_{\infty}$, \[ f_{w}=\pi_{w}\circ
     q:\bbP^1\rightarrow\bbP^k\rightarrow\bbP^1\] where $\pi_{w}$ is
     projection onto a unique line $L_w\subset\bbP^k$ that contains
     $q(w)$.
\end{thm}
We have been unable to improve the theorem from an existence result to
a more satisfying version that would specify $L_w$, and a scale
(described in Section~\ref{sec:rat}), in terms of $q$.

We prove Theorems~\ref{th:sphere}, \ref{th:git} and \ref{th:rat} in
Sections~\ref{sec:sphere}, \ref{sec:git} and \ref{sec:rat}
respectively.  Property (i) of Theorem~\ref{th:sphere} is not
sufficient to guarantee that a curve is the spectral curve of a
monopole.  In general a spectral curve requires even further
restrictions.  Nevertheless, in Section~\ref{sec:ch2} we exploit the
fact that for charge two monopoles no further restrictions are
necessary.  In Section~\ref{sec:massless} we discuss similarities
between the holomorphic sphere and previous work on massless
monopoles. In the final section we prove a vanishing theorem for
hyperbolic monopoles which we need to prove Theorem~\ref{th:sphere}.
This is of some independent interest having been conjectured in
\cite{MSiCom} and is a necessary step towards generalising that work
to monopoles of non-integral mass.

\section{Holomorphic sphere}  \label{sec:sphere}
Theorem~\ref{th:sphere} consists of two quite independent results.
Part (i) states that the spectral curve of a hyperbolic monopole is of
a specific type inside the variety of $(k,k)$ curves in
$\bbP^1\times\bbP^1$.  Part (ii) is a consequence of a more direct
relationship between the spectral curve and the boundary value of the
hyperbolic monopole, given in Theorem~\ref{th:psiherm} in terms of the
defining polynomial of the spectral curve.

The spectral curve of a hyperbolic monopole possesses a type of
positivity property which can be seen explicitly in the case of charge
1 monopoles.  The spectral curve of a charge 1 monopole is a real
$(1,1)$ curve corresponding to all geodesics containing a given point
of $\bbH^3$.  Such a $(1,1)$ curve necessarily lies in the connected
component of the diagonal of $\bbP^1\times\bbP^1$ which represents all
geodesics containing $0\in\bbH^3$.  For example, if a real $(1,1)$
curve contains the points $(0,0)$ and $(\infty,\infty)$, it is of the
form $w-az=0$ for $a\in\bbR^*$.  It is a spectral curve of a charge 1
monopole precisely when $a>0$.  The proof of the first part of
Theorem~\ref{th:sphere} is a generalisation of this simple fact, using
the connectivity of the moduli space and a rather deep analogue of the
property $a\in\bbR^*$.

The defining ``polynomial'' of the spectral curve of a hyperbolic
monopole is an example of a general feature used in this paper.  A
section of
\[ \begin{array}{c}\cO(k,k)\\\downarrow\\Z=\bbP^1\times\bbP^1-\bar{\Delta}
\end{array}\]
extends to a section of 
\[ \begin{array}{c}\cO(k,k)\\\downarrow\\Q=\bbP^1\times\bbP^1
\end{array}\]
and hence is given by a polynomial.  More generally, if a holomorphic
bundle over $Z$ extends to $Q$ then any section extends.  When the
bundle is trivial, this says that a holomorphic function on $Z$ is
necessarily constant, which uses the fact that $Z$ contains many
compact holomorphic curves, in particular those $(1,1)$ curves
corresponding to all geodesics containing a given point of $\bbH^3$.
The more general fact can be proven in a couple of ways.  In the proof
of Lemma~\ref{th:ext}, it is shown that sections of line bundles over
$Z$ lift to homogeneous functions defined over a large enough (to
contain many compact holomorphic curves) subset of $\bbP^3$.  An
alternative argument uses the fact that any local holomorphic function
defined over a deleted neighbourhood of a totally real submanifold, in
this case any open set in $\bar{\Delta}$, extends uniquely to a
holomorphic function on the neighbourhood.  Thus, if a holomorphic
bundle over $Z$ extends to $Q$, then any local holomorphic section
also extends.

\subsection{Positive definite}
Along any geodesic $\gamma\subset\bbH^3$ the monopole $(A,\Phi)$
defines the scattering equation
\begin{equation}  \label{eq:scat}
(\partial_t^A-i\Phi)s=0
\end{equation}
where $t$ parametrises $\gamma$ and $s(t)$ is a section of $E$
restricted to $\gamma$.  The Bogomolny equations define an
integrability condition $[\partial_t^A-i\Phi,\partial_{\bar{z}}^A]=0$
and hence local solutions satisfying $\partial_{\bar{z}}^As=0$ can be
found.  These define a holomorphic bundle $\tilde{E}$ over $Z$ with
distinguished sub-line bundles $L_+,L_-$ given by those solutions that
decay as $t\rightarrow\infty$, respectively $t\rightarrow-\infty$.
The line bundles $L_+$ and $L_-$ coincide over an algebraic curve
$\Sigma\in\bbP^1\times\bbP^1$ known as the {\em spectral curve}.
Points on the spectral curve represent geodesics that possess a
solution which decays both as $t\rightarrow\pm\infty$.

Corresponding to reversing the direction of a geodesic, the space of
geodesics $\bbP^1\times\bbP^1-\bar{\Delta}$ possesses a real structure
$(w,z)\mapsto(\hat{z},\hat{w})$.  The spectral curve
$\Sigma$ is invariant under the real structure since a solution of
(\ref{eq:scat}) along $\gamma(t)$ can be used to construct a solution
of (\ref{eq:scat}) along $\gamma(-t)$ with decay preserved.

\begin{lemma}   \label{th:real}
     The defining polynomial for the spectral curve can be chosen to
     satisfy
     \begin{equation} \label{eq:real}
     \psi(\hat{z},\hat{w})=\overline{\psi(w,z)}
     \end{equation}
     and to be positive on the anti-diagonal $w=\hat{z}$.
\end{lemma}
\begin{proof}
     If we take $\psi(w,z)$ to mean a degree $(k,k)$ polynomial in
     $w^{-1}$ and $z$ (perhaps one would prefer $P(w^{-1},z)$ or to
     refer to $w^k\psi(w,z)$ as the polynomial) then one can express
     the reality condition quite simply.  The reality condition means
     that $\psi(\hat{z},\hat{w})$ and $\psi(w,z)$ have the same zero
     set and since $\psi(w,z)$ and the complex conjugate of
     $\psi(\hat{z},\hat{w})$ both define degree $(k,k)$ polynomials in
     $w^{-1}$ and $z$, they are the same up to a constant
     $\psi(\hat{z},\hat{w})=c\overline{\psi(w,z)}$.  The spectral
     curve does not intersect the anti-diagonal so $\psi$ does not
     vanish there, and hence $c=\exp(2i\theta)$ for some constant
     $\theta$.  We can replace $\psi$ by $\exp(-i\theta)\psi$ to get
     (\ref{eq:real}).  Since $\psi$ does not vanish on the
     anti-diagonal, it is either positive or negative there, and if
     the latter we can replace it by $-\psi$.
\end{proof}

\begin{thm}  \label{th:orth}
     For each monopole $(A,\Phi)$ there exists a holomorphic embedding
     \[ q:\bbP^1\rightarrow\bbP^k\] unique up to the action of
     $U(k+1)$ on its image satisfying \[\langle
     q(\hat{w}),q(z)\rangle=\psi(w,z).\]
\end{thm}
\begin{proof}
     Let $v(z)=(1,z,z^2,\dots,z^k)$.  Then
     \[\psi(w,z)=v(-1/w)^T\Psi v(z)\] for a
     $(k+1)\times(k+1)$ matrix $\Psi$.  Condition (\ref{eq:real}) is
     equivalent to $\Psi=\overline{\Psi}^T$.

     To prove the theorem we will show that the matrix $\Psi$ is
     positive definite so $\Psi=\bar{Q}^TQ$ for an invertible
     $(k+1)\times(k+1)$ matrix $Q$ unique up to $Q\mapsto uQ$ for $u\in
     U(k+1)$.  Then set $q(z)=Qv(z)$, a degree $k$ holomorphic
     map.

     This proves a stronger property of $q$ than simply being an
     embedding---the image of $q$ spans all of $\bbP^k$, and we call it
     {\em full}.  A full map is an embedding since any singular point
     $z$ would satisfy $0=q'(z)=Qv'(z)$, and any double
     point would satisfy $0=q(z_1)-\lambda
     q(z_2)=Q(v(z_1)-\lambda v(z_2))$, and in both cases
     $Q$ would have a non-trivial kernel, contradicting the fullness of
     $q$.

     The following lemma gives part of the property that the bilinear
     form $\Psi$ is positive definite.
\begin{lemma}   \label{th:nondeg}
     The matrix $\Psi$ is non-degenerate.
\end{lemma}
\begin{proof}
     We will prove that for any $(k,k)$ curve $\Sigma=\{(w,z)\in
     \bbP^1\times\bbP^1|\psi(w,z)=0\}$ with coefficient matrix $\Psi$,
     the condition that $\Psi$ be non-degenerate is equivalent to
     $H^0(\Sigma,\cO(k,-2))=0$.  The spectral curve of a mass $m$
     monopole satisfies the property
     $L^{2m+k}|_{\Sigma}\cong\cO_{\Sigma}$ \cite{MSiSpe} so
     \[H^0(\Sigma,\cO(k,-2))=H^0(\Sigma,L^{2m}(0,k-2))\] and the
     latter vanishes by Theorem~\ref{thm:vanishing}.

     A section of $H^0(\Sigma,\cO(k,-2))$ is represented by a
     polynomial together with the defining polynomial of $\Sigma$,
     expressed as $p(z),\psi(z)\in\bbC[w][z]$ with coefficients given
     by sections of $\cO(k)$, or degree $k$ polynomials in $w$, such
     that 
     \begin{equation} \label{eq:section}
     z^2p(z)+\psi(z)q(z)\in\bbC[w][z^{-1}] \quad\mbox{for some}\quad
     q(z)\in\bbC[z,z^{-1}].  
     \end{equation} 
     Put $q=\sum q_lz^{1-l}$.  Then \begin{equation} \label{eq:kernel}
     \psi(w,z)q(z)=\sum_{i,j,l}\Psi_{ij}q_lw^iz^{1+j-l} \end{equation}
     and the coefficient of $w^iz$ is $\sum_j\Psi_{ij}q_j$.  The
     degeneracy of $\Psi$ is equivalent to the existence of a
     non-trivial $q$ such that $\sum_j\Psi_{ij}q_j=0$ for all $i$.
     But then (\ref{eq:kernel}) becomes (\ref{eq:section}) if we move
     the terms on the right hand side of (\ref{eq:kernel}) with
     positive powers of $z$ to the left hand side, and the lemma is
     proven.

     The difference of sections in two charts giving rise to a vector
     in the kernel of $\Psi$ looks like the coboundary map in
     cohomology.  In fact, we can express the proof of the lemma in
     terms of the exact sequence in cohomology given by
     \begin{equation*}
\label{eq:Psi}
	0\rightarrow H^0(\Sigma,\cO(k,-2))\rightarrow
	H^1(Q,\cO(0,-2-k))\stackrel{\Psi}{\rightarrow} H^1(Q,\cO(k,-2))
     \end{equation*}
     where the right-most map is multiplication by $\psi(w,z)$
     and becomes the matrix $\Psi:\bbC^{k+1}\rightarrow\bbC^{k+1}$.
     \end{proof}

     Since a continuous family of non-degenerate Hermitian matrices
     has constant signature, it follows from Lemma~\ref{th:nondeg} and
     the connectivity of the moduli space that we need show that only
     one monopole possesses a positive definite $\Psi$.  This is true
     for axially symmetric monopoles by the explicit construction
     given in Section~\ref{sec:axial} (or we can prove it for
     half-integer mass monopoles using techniques from \cite{ABrBou}.)
     Hence $\Psi$ is positive definite for all monopoles and the
     theorem is proven.
\end{proof}

\subsection{Hermitian metrics}
A Hermitian metric on a vector space $V$ is a linear map 
\[ H:\overline{V}\otimes V\rightarrow\bbC\quad\mbox{satisfying}\quad 
H(u,v)=\overline{H(\bar{v},\bar{u})}\] where the map \large
$\overline{\cdot}$ \normalsize gives an antilinear isomorphism from
$V$ to $\overline{V}$ and back. 

A Hermitian metric on a holomorphic bundle uniquely determines a
Hermitian connection on the holomorphic bundle compatible with the
holomorphic structure.  The reducible connection on the sphere at
infinity $A_{\infty}$, a $U(1)$ connection on the holomorphic line
bundle $\cO(-k)$ over $S^2$, can be described via a Hermitian metric
on $\cO(-k)$.  In local coordinates, the Hermitian metric, $h$, is
locally a positive valued function well-defined up to $h(z)\sim
|g(z)|^2h(z)$, for $g$ a local holomorphic function.  The
$h$-Hermitian connection is $\partial_z\ln h\cdot dz$, or in a unitary
gauge it is
\begin{equation}  \label{eq:metcon} 
A_{\infty}=-\partial_{\bar{z}}\ln\xi\cdot d\bar{z}+\partial_z\ln\xi\cdot dz
\end{equation}
for $\xi^2=h$, the positive square root.

\begin{thm}   \label{th:psiherm}
     Let $\psi$ be the defining polynomial of the spectral curve of
     $(A,\Phi)$, \[\Sigma=\{(w,z)|\psi(w,z)=0\}.\] Then the
     restriction of $\psi$ to the anti-diagonal,
     $\psi|_{\bar{\Delta}}$, gives rise to a Hermitian metric on the
     holomorphic bundle $\cO(-k)$ over $\overline{\Delta}$ that
     defines the connection at infinity.
\end{thm}
\begin{proof}
     The real structure on $\bbP^1\times\bbP^1$ given by
     $(w,z)\mapsto(\hat{z},\hat{w})$ lifts to a real structure on the
     bundle $\cO(k,k)$.  This is reflected in Lemma~\ref{th:real}
     where it is proven that local trivialisations for $\cO(k,k)$ can
     be chosen so that the involution on each fibre is simply complex
     conjugation, and the real structure fixes any section of
     $\cO(k,k)$ whose zero set is preserved by the real structure.

     Any section $s$ of $\cO(k,k)$ gives a map
     $s:\cO(-k,-k)\rightarrow\bbC$.  Suppose $s$ is fixed by the real
     structure.  Restrict $s$ to the fixed point set of the real
     structure, $\overline{\Delta}\subset\bbP^1\times\bbP^1$.  We can
     identify
     $\cO(-k,-k)|_{\overline{\Delta}}\cong\overline{\cO(-k)}\otimes\cO(-k)$
     so $s$ defines a Hermitian metric
     $s:\overline{\cO(-k)}\otimes\cO(-k)\rightarrow\bbC$ on the
     holomorphic bundle $\cO(-k)$ over $\overline{\Delta}$.

     Apply this to $\psi$, the defining polynomial of the spectral
     curve, since it is fixed under the real structure.  Its
     restriction to the anti-diagonal defines a Hermitian metric on
     the holomorphic bundle $\cO(-k)$, and hence a Hermitian
     connection there.  It remains to show that this Hermitian
     connection is the $U(1)$ connection at infinity of the monopole.
     This is a consequence of the following three lemmas.

     In order to understand the map $\psi:\cO(-k,-k)\rightarrow\bbC$
     we choose local holomorphic sections
     $\overline{s_+(\hat{z},\hat{w})}\otimes s_+(w,z)$ of $\cO(-k,-k)$
     where $s_+(w,z)$ is a solution of (\ref{eq:scat}) along the
     geodesic traveling from $\hat{w}$ to $z$ (so
     $s_+(\hat{z},\hat{w})$ is a solution of (\ref{eq:scat}) along the
     oppositely oriented geodesic.)

\begin{lemma}  \label{th:psiact}
     The section $\psi$ acts on $\cO(-k,-k)$ by
     \[\psi(\overline{s_+(\hat{z},\hat{w})}\otimes s_+(w,z))=\langle
     s_+(\hat{z},\hat{w}),s_+(w,z)\rangle.\]
\end{lemma}
\begin{proof}
     Recall from \cite{MSiSpe} that a hyperbolic monopole defines a
     holomorphic bundle $\widetilde{E}\to Z$ with two extensions:
     $$0 \to L^m(0, -k) \to\widetilde{E}\to L^{-m}(0, k) \to 0$$
     and
     $$0 \to L^{-m}(-k,0) \to\widetilde{E}\to L^m(k,0) \to 0$$
     for $L=\cO(1,-1)$.  The sub-bundles $L^+=L^m(0, -k)$ and
     $L^-=L^{-m}(-k,0)$, defined as the space of solutions of
     (\ref{eq:scat}) that decay as $t\rightarrow\infty$, respectively
     $t\rightarrow-\infty$, coincide over the spectral curve $\Sigma$
     and their coincidence defines a non-vanishing section over
     $\Sigma$ of $L^{2m+k}$.

     The spectral curve is a $(k,k)$ curve with defining polynomial
     $\psi$, hence
     \[0\rightarrow\cO(-k,-k)\stackrel{\psi}{\rightarrow}
     \cO\rightarrow\cO_{\Sigma}\rightarrow 0\] we can tensor this
     with $L^m(k,0)$ to get
     \[ \begin{array}{ccccccccc}
     0&\rightarrow&L^m(0,-k)&\stackrel{\psi}\rightarrow
     &L^m(k,0)&\rightarrow&\cO_{\Sigma}(L^m(k,0))&\rightarrow& 0\\
     &&\uparrow\cong&&\uparrow&&\uparrow&&\\
     0&\rightarrow&L^m(0,-k)&\rightarrow&\widetilde{E}
     &\rightarrow&L^{-m}(0,k)&\rightarrow& 0 \end{array}\] 
     which represents $\psi$ as a map
     $\psi:L^+\rightarrow\widetilde{E}/L^-$.

     For $\tau$ the real structure on the space of geodesics, there is
     a natural map \[\widetilde{E}_{\tau(w,z)}\otimes
     \widetilde{E}_{(w,z)}\rightarrow\bbC\] given by $r\otimes
     s\mapsto \langle r(t),s(t)\rangle$ where
     $r\in\widetilde{E}_{\tau(w,z)}$ and $s\in\widetilde{E}_{(w,z)}$
     or equivalently, $(\partial_t^A+i\Phi)r=0$ and
     $(\partial_t^A-i\Phi)s=0$.  The inner product $\langle
     r(t),s(t)\rangle$ is independent of $t$, since
     \[\partial_t\langle r(t),s(t)\rangle=
     \langle(\partial_t^A+i\Phi)r(t), s(t)\rangle+\langle
     r(t),(\partial_t^A-i\Phi)s(t)\rangle=0.\]

     Thus, $\psi:L^+\rightarrow\widetilde{E}/L^-$ can be re-expressed
     as $\psi:(\widetilde{E}/L^-)^*\otimes L^+\rightarrow\bbC$ and
     since $(\widetilde{E}/L^-)^*\cong
     (L^-)^{\bot}=\overline{\tau^*L^+}$, we have
     \[\psi(\overline{s_+(\hat{z},\hat{w})}\otimes s_+(w,z))=\langle
     s_+(\hat{z},\hat{w}),s_+(w,z)\rangle.\]
\end{proof}

\begin{lemma}  \label{th:pert}
     \[\lim_{w\rightarrow\hat{z}}\psi(w,z)=
     \lim_{w\rightarrow\hat{z}}\langle s_+(\hat{z},\hat{w}),s_+(w,z))
     \rangle=\lim_{t\rightarrow\infty}\exp(2mt)\|s_+(w,z)\|^2.\]
\end{lemma}
\begin{proof}
     In order to make sense of the lemma, we really need to choose
     local trivialisations for the bundles so that we are dealing with
     local functions, and so that $s_+$ is well-defined.  We can do
     this as follows.  Fix $\lim_{t\rightarrow\infty}\exp(mt)s_+(w,z)$
     in a small neighbourhood $U\subset\bbP^1$ containing $z$.  (One
     can choose any family of solutions $s_+(w_0,z)$, for fixed $w_0$,
     and use $\lim_{t\rightarrow\infty}\exp(mt)s_+(w_0,z)$.)  As $w$
     moves close enough to $\hat{z}$ so that $\hat{w}\in U$, we use
     the same limit for
     $\lim_{t\rightarrow-\infty}\exp(mt)s_+(\hat{z},\hat{w})$, where
     we use the parameter $-t$ for the oppositely oriented geodesic.
     
     The proof is not yet immediate, since we have only arranged that
     the values at opposite ends of a geodesic in the $t$ independent
     quantity $\langle s_+(\hat{z},\hat{w}),s_+(w,z))$ are
     approximately the same.

     Now use the fact from \cite{RadSin} that there exists a gauge in
     which \begin{equation} \label{eq:pert} \partial_t^A\pm
     i\Phi=\partial_t\pm i\left(\begin{array}{cc}im&0\\0&-im
     \end{array}\right)+\epsilon\cdot C\exp(-m|t|) \end{equation}
     where $C$ is constant and $\epsilon\rightarrow 0$ as
     $w\rightarrow z$.  We see that we do indeed end up with the
     product of
     $\lim_{t\rightarrow-\infty}\exp(mt)s_+(\hat{z},\hat{w})$ and
     $\lim_{t\rightarrow\infty}\exp(mt)s_+(w,z)$ which, by
     construction, tends towards
     $\lim_{t\rightarrow\infty}\exp(2mt)\|s_+(w,z)\|^2$.  Further
     details can be found in \cite{NorBou}.
\end{proof}

\begin{lemma}  \label{th:ainf}
     Let $s_+$ be a local holomorphic section of $L_+$.  Then if we
     fix $w$ and parametrise the sphere at infinity by $z$
     \[h(w,z)=\lim_{t\rightarrow\infty}\exp(2mt)\| s_+(w,z)\|^2\] is a
     Hermitian metric that determines the $U(1)$ connection at
     infinity.
\end{lemma}
\begin{proof}
     Since the $U(1)$ connection at infinity is Hermitian, it can be
     determined from its $(0,1)$ part.  If we fix one end of a family
     of geodesics (to be $\hat{w}$) and vary the other end ($z$), then
     $\partial_{\bar{z}}^As_+=\lambda(z)s_+$ for some $\lambda(z)$
     {\em independent of} $t$.  This follows from the three properties
     $(\partial_t^A-i\Phi)s_+=0$,
     $[\partial_{\bar{z}}^A,\partial_t^A-i\Phi]=0$ and
     $\partial_{\bar{z}}^As_+$ decays as $t\rightarrow\infty$.  In
     particular, $\lambda(z)$ makes sense at $t=\infty$ and gives the
     $(0,1)$ part of the $U(1)$ connection at infinity.  The $(1,0)$
     part can be determined by the fact that the $U(1)$ connection at
     infinity is Hermitian with respect to
     $\lim_{t\rightarrow\infty}\exp(2mt)\| s_+(w,z)\|^2$.
\end{proof}
     From Lemmas~\ref{th:psiact}, \ref{th:pert} and \ref{th:ainf}, we
     see that $\psi|_{\bar{\Delta}}$ defines a Hermitian metric that
     gives rise to the $U(1)$ connection at infinity.

     Note that when choosing a local frame for $\cO(-k,-k)$, if we
     also require
     \[\partial_{\bar{z}}^As_+(w,z)=0=\partial_w^As_+(w,z)\] and
     similar conditions on $s_+(\hat{z},\hat{w})$, then local
     holomorphic sections for $\cO(-k,-k)$ are simply given in terms
     of local holomorphic functions with respect to this frame,
     whereas without these extra conditions one must use the $(0,1)$
     part of a connection to detect local holomorphic sections.

     It is important to understand that in order to use
     (\ref{eq:metcon}) to retrieve the connection from the Hermitian
     metric, one needs a local trivialisation of the holomorphic
     bundle in which local holomorphic sections are given by local
     holomorphic functions.  When we choose {\em separable} transition
     functions for $\cO(-k,-k)$, that is each transition function is
     given by a product of transition functions for $\cO(-k,0)$ and
     $\cO(0,-k)$, then the holomorphic structure on the bundle
     $\cO(-k)$ over the antidiagonal has local holomorphic sections
     given by local holomorphic functions.  In particular, the choice
     of $\psi$ as a {\em polynomial} (in the local coordinates $w$ and
     $z$ or $-1/w$ and $z$, etc) arises from separable transition
     functions.  Thus, we can choose $\psi$ to be the defining
     polynomial and use (\ref{eq:metcon}) to retrieve the connection
     and Theorem~\ref{th:psiherm} is proven.
\end{proof}

In the statement of Theorem~\ref{th:sphere}, we express the
relationship of the holomorphic map $q:\bbP^1\rightarrow\bbP^k$ with
the connection at infinity via $F_{A_{\infty}}=q^*\omega$, for
$\omega$ the Kahler form on $\bbP^k$.  The holomorphic map $q$ pulls
back a Hermitian metric, its connection and its curvature.  The
Hermitian metric is given by $\langle q(z),q(z)\rangle$ which, by
Theorem~\ref{th:orth} is $\psi|_{\bar{\Delta}}$.  Thus,
Theorem~\ref{th:psiherm} can be restated as $q$ pulls back the
Hermitian metric that defines the $U(1)$ connection at infinity and in
particular $F_{A_{\infty}}=q^*\omega$ and part (ii) of
Theorem~\ref{th:sphere} is proven.

The real analyticity of $\psi$ in a neighbourhood of ${\bar{\Delta}}$,
and the fact that ${\bar{\Delta}}$ is a totally real submanifold of
$\bbP^1\times\bbP^1$ allows one to show that $\psi|_{\bar{\Delta}}$
well-defined up to multiplication by holomorphic and anti-holomorphic
functions uniquely extends and hence determines $\psi$ on
$\bbP^1\times\bbP^1$.  See \cite{NorBou} for details.  Thus
\begin{cor}  \label{th:infdet}
     An $SU(2)$ hyperbolic monopole is determined up to gauge by its
     reducible connection on the sphere at infinity.
\end{cor}
This was proven in \cite{NorBou} by a slightly different method.  That
paper did not require
Theorem~\ref{th:psiherm}---$\psi|_{\bar{\Delta}}$ is a Hermitian
metric that defines the connection at infinity---although it did use
the fact that $\psi|_{\bar{\Delta}}$ determines $\psi$.  The proof of
Corollary~\ref{th:infdet} in the half integer mass case \cite{ABrBou}
uses the discrete Nahm equations to prove that a holomorphic map
$q:\bbP^1\rightarrow\bbP^k$ uniquely determines the monopole, and a
result of Calabi \cite{CalIso} to show that $q^*\omega$ uniquely
determines $q$.  That approach, combined with Theorem~\ref{th:sphere}
can be used to give a third proof of Corollary~\ref{th:infdet},
although the use of Calabi's theorem seems a bit unnecessary given the
alternative local argument.

\section{Centred monopoles}  \label{sec:git}
The holomorphic sphere $q:\bbP^1\rightarrow\bbP^k$ associated to a
monopole allows one to use geometric invariant theory to define the
centre of a monopole.  We can represent $q$ by
\[ q(z)=v_0+\sqrt{k}v_1z+\sqrt{k\choose 2}v_2z^2+\dots+
\sqrt{k\choose j}v_jz^j+\dots +v_kz^k\]
where each $v_j\in\bbC^{k+1}$ and the coefficients arise quite
naturally as we shall see later.  The $PSL(2,\bbC)$ action on the
domain of such maps lifts to a linear action of $SL(2,\bbC)$ on
$(k+1)$-tuples $(v_0,\dots,v_k)\mapsto(w_0,\dots,w_k)$ where,
\begin{eqnarray*}
\left(\begin{array}{cc}a&b\\c&d\end{array}\right)\cdot q(z)&=&
(cz+d)^kv_0+\dots+(cz+d)^{k-j}(az+b)^j\sqrt{k\choose j}v_j+\dots\\
&=&w_0+\sqrt{k}w_1z+\dots+\sqrt{k\choose j}w_jz^j+\dots +w_kz^k.
\end{eqnarray*}
The space of $(k+1)$-tuples is a subset of $\bbC^N$ (for $N=(k+1)^2$),
so that geometric invariant applies.  The norm on the space of
$(k+1)$-tuples is
\begin{equation}  \label{eq:norm}
\|q\|^2=\|v_0\|^2+\|v_1\|^2+\dots+\|v_j\|^2+\dots+\|v_k\|^2
\end{equation}
which is preserved by $SU(2)\subset SL(2,\bbC)$.  (We are abusing
notation by labeling $(v_0,\dots,v_k)\in\bbC^N$ by $q$ when really
$q$ is the projective class in $\bbC\bbP^{N-1}$.)

Recall that a $(k+1)-$tuple $(v_0,\dots,v_k)$ is {\em stable} under
the $SL(2,\bbC)$ action if and only if the map
$SL(2,\bbC)\rightarrow\bbC^N$ given by $g\mapsto
g\cdot(v_0,\dots,v_k)$ is proper, so in particular the $SL(2,\bbC)$
orbit is closed and we can minimise the norm $\|q\|$ in its
$SL(2,\bbC)$ orbit.

\begin{lemma}   \label{th:stable}
Each $(k+1)-$tuple $(v_0,\dots,v_k)$ arising from a degree $k$
holomorphic map $q$ is a stable point of the $SL(2,\bbC)$ action.
\end{lemma}
\begin{proof}
By the Hilbert criterion it is enough to test the stability of a point
on one-parameter subgroups of $SL(2,\bbC)$.  Any one parameter
subgroup in $SL(2,\bbC)$ is given by
\[ g\left(\begin{array}{cc}t&0\\0&t^{-1}\end{array}\right)g^{-1}.\]
Since the degree of $q$ is $k$ then $v_0\neq 0$ and $v_k\neq 0$.
Also, after acting by $g$ the map $q$ is still of degree $k$ and hence
we may assume that $v_0\neq 0$, $v_k\neq 0$ and $g=I$.  Then the
action is given by
\[ (v_0,\dots,v_k)\mapsto(t^{-k}v_0,\dots,t^{2j-k}v_j,\dots,t^kv_k).\]
In particular, since $v_0\neq 0$ and $v_k\neq 0$, the norm
$\|q\|\rightarrow\infty$ as $t\rightarrow 0$ and $t\rightarrow\infty$
so the map is proper.
\end{proof}
\begin{prop}  \label{th:moment}
The moment map for the action of $SU(2)$ is
\[\mu(v_0,\dots,v_k)=(\sum_{j=0}^k(2j-k)\|v_j\|^2,\sum_{j=0}^{k-1}
\sqrt{(j+1)(k-j)}( v_j,v_{j+1}))\in\bbR\times\bbC.\]
\end{prop}
\begin{proof}
We wish to minimise the norm (\ref{eq:norm}) on each $SL(2,\bbC)$
orbit (which is closed by Lemma~\ref{th:stable}.)  The minimum occurs
on stationary points of the infinitesimal action of ${\bf
sl}(2,\bbC)$.  Since ${\bf su}(2)\subset{\bf sl}(2,\bbC)$ acts
trivially it is enough to consider the actions of
\[e_0=\left(\begin{array}{cc}1&0\\0&-1\end{array}\right),\
e_{\theta}=\left(\begin{array}{cc}0&\exp(i\theta)\\0&0\end{array}\right)\]
given by
\[ e_0\cdot(v_0,\dots,v_k)=(-kv_0,(2-k)v_1,\dots,(2j-k)v_j,\dots,kv_k)\]
\[e_{\theta}\cdot(v_0,\dots,v_k)=\exp(i\theta)
(\sqrt{k}v_1,\dots,\sqrt{j(k+1-j)}v_j,\dots,\sqrt{k}v_k,0).\]
Then
\[ e_0\cdot\|q\|^2=2\sum_{j=0}^k(2j-k)\|v_j\|^2\]
\[ e_{\theta}\cdot\|q\|^2=2Re\exp(i\theta)\sum_{j=0}^{k-1}
\sqrt{(j+1)(k-j)}( v_j,v_{j+1})\]
and the result follows.
\end{proof}

\begin{defn}
An $SU(2)$ hyperbolic monopole is {\em centred} at $0\in\bbH^3$ if its
associated holomorphic sphere $q:\bbP^1\rightarrow\bbP^k$ lies in the
zero set of the moment map $\mu$.
\end{defn}

A consequence of the preceeding definition and the discussion of
geometric invariant theory is a well-defined centre of a monopole.
Each $PSL(2,\bbC)$ orbit of a monopole possesses a unique $SO(3)$
orbit that lies in the zero set of the moment map $\mu$.  Hence to
each element in a $PSL(2,\bbC)$ orbit one can associate a unique point
of $\bbH^3$ which is defined to be the centre of the monopole.

\section{Rational maps}  \label{sec:rat}
{\em Proof of Theorem~\ref{th:rat}}.
The rational map $f_{w}(z):\bbP^1\rightarrow\bbP^1$ is defined as
follows.  Consider all geodesics that begin at $w\in S^2_{\infty}$.
Frame the bundle $E$ at $w\in S^2_{\infty}$ and extend it to a
neighbourhood.  This consists of choosing vectors in each of the
eigenspaces of $\Phi$.  One of these vectors extends to a unique
global solution of (\ref{eq:scat}) as follows.  Define $s_{w}(z)$ to
be a solution of (\ref{eq:scat}) along all geodesics beginning at
$w\in S^2_{\infty}$ and extending to $S^2$, satisfying

(i) $\lim_{t\rightarrow-\infty}e^{mt}s_w(z)$ is a non-zero
vector in the chosen eigenspace of $\Phi$;

(ii) $\lim_{z}\rightarrow ws_w(z)$ exists;

(iii) $\partial_{\bar{z}}^As_w(z)=0$.

Similarly, the other eigenspace gives rise to decaying solutions $s_-$
that satisfy these conditions with $e^{mt}$ in (i) replaced by
$e^{-mt}$.  This frame of solutions is unique since any other frame
differs by a holomorphic gauge transformation defined over $S^2$ and
hence is constant, and in fact the identity since the bundle $E$ is
framed at $w$.

Amongst solutions of (\ref{eq:scat}) along each geodesic that begins
at $w\in S^2_{\infty}$ is a solution $s_+$ that decays so
that $\lim_{t\rightarrow\infty}e^{mt}s_+$ is well-defined.  This
defines a one-dimensional subspace of the frame defined above, and
hence of $\bbC^2$.

Thus, we get a map $f_w(z):\bbP^1\rightarrow\bbP^1$ which turns out
to be holomorphic \cite{AtiIns,AtiMag}.  The poles of $f_{w}(z)$
correspond to those points $z_i$ such that the solution $s_+$ along
the geodesic from $w$ to $z$ decays at both ends.  Equivalently, $s_+$
is a multiple of $s_-$ and has no $s_w$ component.  Thus, the poles
come from points of the spectral curve, $(\hat{w},z_i)\in\Sigma$.  We
have chosen a direction in the frame $\bbC^2$ to represent $\infty$.
We choose the orthogonal direction in $\bbC^2$ to represent $0$, so a
zero of $f_{w}$ corresponds to a solution $s_+$ that is a multiple of
$s_w$ and thus has no $s_-$ component.  In particular, $w$ is a zero
of $f_w$ since in the limit $z\rightarrow w$, $s_+$ and $s_-$ are
orthogonal.

Choose $w\in S^2$ and let $L$ be any line in $\bbP^k$ that contains
the point $q(w)$.  Let $P:\bbP^k\rightarrow L$ be projection onto the
line.  It is alternatively described as projection onto the plane in
$\bbC^{k+1}$ defined by $L$ using the Hermitian product on
$\bbC^{k+1}$.  The map $Pq(z):\bbP^1\rightarrow\bbP^1$ is a degree $k$
holomorphic map.  We choose the direction $q(w)\in L$ to represent the
point $0$, and the orthogonal direction to represent $\infty$.  Thus
$Pq(z)$ has poles given by $z_i$ such that $Pq(z_i)\in L$ is
orthogonal to $q(w)\in L$, so $\langle q(w),q(z_i)\rangle=0$.  The
poles correspond to points of the spectral curve
$(\hat{w},z_i)\in\Sigma$ and coincide with the poles of the rational
map.

Furthermore, since $Pq(w)=q(w)$, $w$ is a zero of $Pq(z)$ which agrees
with $f_w(w)=0$.  There are $k$ zeros $\{ w_i|i=1,..,k\}$ ($w_1=w$) of
$f_w$ counted with multiplicity.  When the zeros are distinct,
$q(w_i)$ define a $k$-dimensional subspace of $\bbC^{k+1}$, since $q$
is full, and this possesses a unique orthogonal direction.  Choose
$L_w$ to represent the plane spanned by this orthogonal direction and
$q(w)$.  Hence, the holomorphic map $Pq(z)$ has the same zeros and
poles as $f_w(z)$.  If a zero $w_i$ has multiplicity $d+1$, then
$q(w_i),q'(w_i),..,q^{(d)}(w_i)$ spans a $(d+1)-$dimensional subspace
of $\bbC^{k+1}$ and the unique orthogonal direction still exists.

The rational maps $f_w(z)$ and $Pq(w)$ differ by a constant.  This
constant determines a scale on the line $L_w$.  We supposed that the
coefficients of the unit vector in $q(w)$ and the orthogonal unit
vector in $L_w$ determine a rational map, or in other words that the
isomorphism of $L_w$ with $\bbP^1$ respects the metric on $\bbP^k$.
It may be that there is another natural scale on $L_w$.  The
question of how we might determine $L_w$ and the scale intrinsically
from $q$ is an interesting one.

\section{Charge two monopoles}   \label{sec:ch2}
For charge two hyperbolic monopoles, we can get explicit expressions
for the boundary data.  We will restrict to the space of centred
charge 2 hyperbolic monopoles since these give rise to interesting
structure.  A charge two monopole is centred if after reflection
in the origin, the new monopole is gauge equivalent to the original
one.  Since a gauge equivalent monopole produces the same equivalence
class of holomorphic spheres in projective space,
$q:\bbP^1\rightarrow\bbP^2$, comes from a centred monopole when
\begin{equation}  \label{eq:holc}
q(\hat{z})=u\cdot q(z),\ \ {\rm for\ some}\ u\in U(3)
\end{equation}
($u$ is independent of $z$.)  If we put $q=v_0+v_1z\sqrt{2}+v_2z^2$ for
$v_i\in\bbC^3$, then (\ref{eq:holc}) is equivalent to
$\|v_0\|^2=\|v_2\|^2$ and $( v_0,v_1)+(
v_1,v_2)=0$ and this is the zero set of the moment map defined in
Proposition~\ref{th:moment}.  Put
\begin{equation}  \label{eq:m2}
\M_2=\{ (v_0,v_1,v_2)\in\bbC^3\otimes\bbC^3\ |\
\|v_0\|^2=\|v_2\|^2,\ ( v_0,v_1)+( v_1,v_2)=0\}/CU(3)
\end{equation}
where $CU(3)=\bbR^+\times U(3)$ is the conformal unitary group which
acts on a triple by $(v_0,v_1,v_2)\mapsto (u^{-1}v_0,u^{-1}v_1,u^{-1}v_2)$.

If we replace the vectors $v_i$ in (\ref{eq:m2}) by vectors in
$\bbC^2$ and quotient by $CU(2)$ then this gives the space of centred
rational maps of degree 2 which naturally sit inside $\M_2$.  if we
replace the vectors $v_i$ in (\ref{eq:m2}) by vectors in $\bbC$ and
quotient by $\bbC^*$ then we get a real structure on $\bbC\bbP^2$ with
ficed set $\bbR\bbP^2$.

The space $\M_2$ is a five dimensional space that contains an open
dense five-dimensional manifold $\M_2^0\subset\M_2$ that is given by
triples of independent vectors.  Points of $\M_2^0$ precisely
correspond to full maps $q$ and these contain the space of centred
charge 2 hyperbolic monopoles.

There is an $SO(3)$ action on $\M_2$ that preserves $\M_2^0$ coming
from the action of $SU(2)$ on the polynomials $(1,z,z^2)$ given by
\[\left(\begin{array}{cc}a&b\\-\bar{b}&\bar{a}\end{array}\right)\cdot
(1,z,z^2)=((-\bar{b}z+\bar{a})^2,(-\bar{b}z+\bar{a})(az+b),(az+b)^2).\]
It is well-defined since it commutes with the $CU(3)$ action.  A
convenient description of the space $\M_2$ is as follows.
\begin{prop}  \label{th:isom}
     $\M_2\cong{\bf su}(2)\otimes{\bf su}(2)/CO(3)$ and the isomorphism
     respects the right $SO(3)$ actions.
\end{prop}
\begin{proof}
     Again $CO(3)=\bbR^+\times SO(3)$ is the conformal orthogonal
     group.  A point of ${\bf su}(2)\otimes{\bf su}(2)/CO(3)$ is
     represented by a triple $(r_0,r_1,r_2)$ for $r_i\in{\bf su}(2)$,
     and the isomorphism is given by
     \[(r_0,r_1,r_2)\mapsto\left(\frac{1}{\sqrt{2}}(r_0+r_2i),r_1,
     \frac{1}{\sqrt{2}}(-r_0+r_2i)\right).\]
     The proof requires the choice of representatives in each $CU(3)$
     orbit \[ \M_2\cong\{(v_0,v_1,-\bar{v}_0)\ |\ v_1\in\bbR^3\}.\]

     Using $CU(3)$, we may assume that $v_1=(1,0,0)$ so that by
     (\ref{eq:m2}) $v_0=(c,\xi_0)$ and $v_2=(-\bar{c},\xi_2)$ for
     $\xi_i\in\bbC^2$ satisfying $\|xi_0\|^2=\|xi_2\|^2$.  Now use
     $u\in U(2)$ to realise $u\xi_2=-\bar{u}\bar{\xi}_0$, or
     equivalently $u^Tu\xi_2=-\bar{\xi}_0$.  We
     can do this since $\{ u^Tu|u\in U(2)\}$ acts transitively on
     $S^3\subset\bbC^2$.
\end{proof}

One would expect that $\M_2^0$ is a one parameter family of
four-dimensional manifolds of centred charge 2 hyperbolic monopoles
with given mass.  In fact, it seems that half of $\M_2^0$ does not
represent hyperbolic monopoles.  Evidence for this is the fact that
the point of $\M_2^0$ consisting of $\{v_i=e_i|i=0,1,2\}$, where
$e_0,e_1,e_2$ is an orthonormal set of basis vectors, is the unique
fixed point of the $SO(3)$ action on $\M_2^0$.  Such a point cannot
correspond to a hyperbolic monopole, since no monopole is $SO(3)$
invariant.  It does correspond to all charge 2 Euclidean monopoles
since they each give a symmetric measure at infinity which would be
pulled back by the this fixed point.

Consider the axially symmetric points in $\M_2^0$.  These are given by
orthogonal triples $(v_0,v_1,v_2)$ and thus the spectral curve is
given by
\[w^2-2\|v_1\|^2wz+z^2=0.\]
In Section~\ref{sec:axial} we calculate the spectral curves of axially
symmetric monopoles.  In the charge 2 case, we find that the spectral
curve is $w^2-2\cos(\pi/(2m+2))wz+z^2$.  Thus
$\|v_1\|^2=\cos(\pi/(2m+2))$ and in particular it takes its values on
the unit interval and the symmetric point is on one side of the
allowed values.

In general, we expect the four dimensional spaces of monopoles with
given mass to form two sided hypersurfaces in $\M_2^0$.  We expect the
symmetric point to partition $\M_2^0$ into two pieces, one containing
hyperbolic monopoles.  The piece containing hyperbolic monopoles is
determined by the axially symmetric examples, and by the fact that in
the limit, as the triple tends towards spanning a two-dimensional
subspace, the massless monopoles emerge.  One might guess that the
other half of the points correspond to asymptotic values of spherical
monopoles near a singular point.

To identify the mass of the monopole from the point of $\M_2^0$ is
difficult.  However, we can in a sense understand the tangential
direction of changing mass as follows.  Associated to a hyperbolic
monopole is the rational map obtained from scattering from
$0\in\bbH^3$ \cite{JNoCom}, and when the monopole has charge 2 and is
centred, it is uniquely determined by the intersection of the spectral
curve with the diagonal in $\bbP^1\times\bbP^1$.  We fix this rational
map, and change the mass.

On ${\bf su}(2)\otimes{\bf su}(2)/CO(3)$ for $\nu=(r_0,r_1,r_2)$ the
map $\nu\mapsto[\nu,\nu]$ is well-defined.  Here $[\cdot,\cdot]$ is
the bracket induced on ${\bf su(2)}\otimes{\bf su(2)}$ by the Lie
bracket on the ${\bf su(2)}$.  Note that this is not a Lie bracket,
and in general $[\nu,\nu]\neq 0$.

\begin{prop}
     The map $\nu\mapsto[\nu,\nu]$ is an involution with fixed point
     the symmetric point of $\M_2^0$.
\end{prop}
\begin{proof}
     Put $\nu=e_1\otimes r_1+e_2\otimes r_2+e_3\otimes r_3$ where
     $e_i$ is an orthonormal basis of ${\bf su}(2)$.  Then
     \[[\nu,\nu]=2e_1\otimes[r_2,r_3]+2e_2\otimes[r_3,r_1]+2e_3\otimes
     [r_1,r_2]\]
     since $[e_i,e_j]=\epsilon_{ijk}e_k$.  The element $\nu$ is fixed
     by this map if $[r_i,r_j]=\lambda\epsilon_{ijk}r_k$ for a
     constant $\lambda$, thus $r_i\mapsto e_i$ under the action of
     $CO(3)$.  The square of this map is given by
     \begin{eqnarray*}
	[[\nu,\nu],[\nu,\nu]]&=&8e_1\otimes[[r_3,r_1],[r_1,r_2]]+8e_2
         \otimes[[r_1,r_2],[r_2,r_3]]\\&&\ \ \ \ \ \ \ \ \ \ \ \ \ \ \
         +8e_3\otimes[[r_2,r_3], [r_3,r_1]]\\
	&=&8\langle r_1,[r_2,r_3]\rangle\nu\\
	&\equiv&\nu
     \end{eqnarray*}
     where the last equivalence uses the fact that for monopoles
     $r_1$, $r_2$ and $r_3$ are independent and hence $\langle r_1,[r_2,
     r_3]\rangle r_1\neq 0$.  We have also used the identity
     \[[[r_1,r_2],[r_2,r_3]]=\langle r_1,[r_2,r_3]\rangle r_1\]
     which can be shown to hold in ${\bf su}(2)$ by linearly extending
     the easy identity \[[[e_i,e_j],[e_j,e_k]]=\epsilon_{ijk}e_j.\]

\end{proof}

\begin{prop}
     The expression $[\nu,\nu]$ defines a vector field on $\M_2^0$ and
     a flow along that vector field fixes the rational map of the
     monopole.
\end{prop}
\begin{proof}
     The intersection of the spectral curve with the diagonal in
     $\bbP^1\times\bbP^1$ consists of four points, given by two pairs
     of antipodal points on the diagonal.  These uniquely determine the
     the rational map obtained by scattering from $0\in\bbH^3$.

     Put $v_0=(1/\sqrt{2})(r_0+r_2i)$, $v_1=r_1$ and
     $v_2=(1/\sqrt{2})(-r_0+r_2i)$ for real vectors $r_i$ as in
     Proposition~\ref{th:isom}.  When we translate the triple
     $\nu=(r_0,r_1,r_2)$ in the $[\nu,\nu]$ direction the intersection
     of the spectral curve with the diagonal is unchanged and hence the
     rational map is preserved.  To see this, calculate $\langle
     q(\hat{z}),q(z)\rangle=0$ in terms of the $r_i$ to get the degree
     4 polynomial
     \begin{eqnarray*}
	0&=&\frac{1}{2}\left[\langle r_2,r_2\rangle-\langle r_0,r_0\rangle
	+2i\langle r_0,r_2\rangle\right]z^4
	+\frac{1}{2}\left[\langle r_2,r_2\rangle-\langle r_0,r_0\rangle
	-2i\langle r_0,r_2\rangle\right]\\
	&&+2\left[\langle r_0,r_1\rangle-i\langle r_2,r_1\rangle\right]z^3
	-2\left[\langle r_0,r_1\rangle+i\langle	r_2,r_1\rangle\right]z\\
	&&+\left[\langle r_0,r_0\rangle +\langle r_2,r_2\rangle-2\langle
	r_1,r_1\rangle\right]z^2.
     \end{eqnarray*}
     Now, consider the infinitesimal change given by
     $\nu\mapsto\nu+t[\nu,\nu]$ so
     \[ r_0\mapsto r_0+t[r_1,r_2],\ r_1\mapsto r_1+t[r_2,r_0],\
     r_2\mapsto r_2+t[r_0,r_1].\]
     Up to first order, this yields the change
     \[\langle r_i,r_j\rangle\mapsto\langle r_i,r_j\rangle
     +t\delta_{ij}\epsilon_{ikl}\langle r_i,[r_k,r_l]\rangle\]
     where we only sum over $k$ and $l$.  Thus the coefficients of the
     degree 4 polynomial are unchanged up to first order.  (For example
     take the coefficient of $z^4$,
     \begin{eqnarray*}
	\langle r_0,r_0\rangle&\mapsto& \langle r_0,r_0\rangle
	+2t\langle r_0,[r_1,r_2]\rangle\\
	\langle r_2,r_2\rangle&\mapsto& \langle r_2,r_2\rangle
	+2t\langle r_2,[r_0,r_1]\rangle\\
	\langle r_0,r_2\rangle&\mapsto& \langle r_0,r_2\rangle
     \end{eqnarray*}
     and the changes cancel.)
\end{proof}

For any $\xi\in{\bf so}(3)$, $[\xi,\nu]$ consists of trivial vectors
(they point in the gauge direction) and $[\nu,\xi]$ gives the tangent
space to the $SO(3)$ action.  Thus we have four tangent directions in
$\M_2^0$, three tangent to a moduli space of monopoles with fixed
mass, and one ``transverse'' to each moduli space.  It would be useful
to find another mass-preserving tangent direction that would enable
one to specify the fixed mass submanifold of $\M_2^0$.

\subsection{Mass of the monopole.}

The spectral curve $\Sigma$ of a charge $k$ monopole is a real $(k,k)$
curve in $\bbP^1\times\bbP^1$ with the extra condition that
$\cO(-(2m+k),2m+k)|\Sigma\cong\cO$, \cite{MSiSpe}, where $m$ is the
mass of the monopole, $m=\lim_{r\rightarrow\infty}\|\Phi\|$.  It is
quite difficult to detect the mass from the spectral curve.

The case of charge 2 $SU(2)$ hyperbolic monopoles is special mainly
because it is related to elliptic functions via its elliptic spectral
curve, and because the spectral curve is identified with its Jacobian,
the place where the Nahm data resides.  Elliptic functions and
isomonodromic deformations are used in \cite{HitNew} to find a new
family of Einstein metrics and explicit expressions for them on the
space of charge 2 centred monopoles, those monopoles invariant (up to
gauge) under reflection in the origin.  

We will describe part of the construction in \cite{HitNew}.  For a
generic choice of $(w,z)\in\Sigma$, the two lines $\{w\}\times\bbP^1$
and $\bbP^1\times\{ z\}$ meet $\Sigma$ again, once each.  Label
$(w,z)$ by $P_0$ and the other intersection point of the vertical line
$\{w\}\times\bbP^1$ with $\Sigma$ by $P_1$.  At $P_1$ take a
horizontal line and label the second point of $\Sigma$ which it
intersects by $P_2$.  Continue this process until $P_{4m+4}$ to get
$P_0,P_1,\dots P_{4m+4}$.  Each point $P_i$ gives a divisor on
$\Sigma$, and $P_0+P_1\sim\cO(0,1)$, $P_1+P_2\sim\cO(1,0)$,
$P_2+P_3\sim\cO(0,1)$ and so on, where $\cO(0,1)$ and $\cO(1,0)$ mean
the restriction of these line bundles to $\Sigma$.  Take the
alternating sum of these divisors to get
$P_0+P_1-(P_1+P_2)+\dots-(P_{4m+3}+P_{4m+4})\sim\cO(2m+2,(-2m+2))$.
But $\cO(2m+2,(-2m+2))\sim\cO$ on $\Sigma$ (where we assume for the
moment that $m\in(1/2)*\bbZ$.)  Thus $P_0-P_{4m+4}$ is the trivial
divisor and hence there is a meromorphic function on $\Sigma$ with
multiplicity one zero and pole given respectively by $P_0$ and
$P_{4m+4}$.  Non-constant meromorphic functions must have at least two
zeros, thus we get a contradiction unless $P_0=P_{4m+4}$.

Out of interest, we will mention the relation of this construction to
the Poncelet polygon problem---to find $n$-sided polygons in the plane
inscribed in one conic and circumscribed about another---described in
\cite{HitNew}.  Consider the map
$\pi:\bbP^1\times\bbP^1\rightarrow\bbP^2$ defined by
$\pi((w_0,w_1),(z_0,z_1))=(w_0z_0,w_0z_1+w_1z_0,w_1z_1)$, (or affinely
$\pi(w,z)=(wz,w+z)$.)  The preimage of any point consists of $(w,z)$
and $(z,w)$ thus the map is a two fold branched cover ramified on the
diagonal and branched over the conic $B=(z_0^2,2z_0z_1,z_1^2)$.  It
simply relates the coefficients of a degree two polynomial to its
roots.  The image of any vertical or horizontal line $\{
w\}\times\bbP^1$, respectively $\bbP^1\times\{ z\}$, is tangent to
$B$.

The spectral curve of a {\em centred} 2 monopole is invariant under
the involution that swaps the two factors $(w,z)\mapsto(z,w)$.  This
is because the two points represent a geodesic running from $\hat{w}$
to $z$, respectively a geodesic running from $\hat{z}$ to $w$.  These
geodesics are images of each other under reflection in the origin.
The image of the spectral curve of a centred 2 monopole is a conic,
$C=\pi(\Sigma)$.  Let $(w,z)\in\Sigma$, then the images of the two
lines $\{w\}\times\bbP^1$ and $\bbP^1\times\{ z\}$ are the two
tangents of the conic $B$ meeting $C$ in the point $\pi(w,z)$.  Hence
the construction described above yields a $4m+4$-sided polygon in the
plane inscribed in one conic and circumscribed about another.

When $m\in(1/4)*\bbZ$ one can conclude that $P_0\neq P_{4m+4}$ but
they lie in the same fibre of $\pi$, again giving a solution of the
Poncelet problem.

This construction can be interpreted in terms of the holomorphic
sphere $q:\bbP^1\rightarrow\bbP^2$.  A generic point $z_0\in\bbP^1$
gives rise to a sequence of points $\dots z_{-1},z_0,z_1,z_2,\dots$ by
requiring the condition that $( q(z_i),q(z_{i\pm 1}))=0$.

\begin{lemma}
When $m\in\bbQ$, then the sequence $\{\dots
z_{-1},z_0,z_1,z_2,\dots\}$ defined by $( q(z_i),q(z_{i\pm 1}))=0$ is
a discrete lattice on the sphere.
\end{lemma}
\begin{proof}
This follows from the argument described above.  If $m=m_1/m_2$ then
after $4m_1+4m_2$ steps we can conclude from
$\cO(2m+2,(-2m+2))^{m_2}\sim\cO$ on $\Sigma$ that the sequence closes up.
\end{proof}

The number of points in the discrete lattice is related to the sum of
the numerator and denominator of the mass.  It would be good to see
the mass precisely from the lattice and to understand what can be done
in the irrational mass case.

A lattice can be constructed in this way for any charge $k$ monopole.
At each new step, $k$ new points on the sphere are produced.  It is
unlikely that this will yield a discrete lattice.  The argument for
this relied crucially on the property that after a finite number of
steps in the construction, we are left with a question about the
divisor $P_0-P_N$ consisting of two points, and can use the fact that
a meromorphic function must have at least two zeros.

Hitchin mentions that his metrics are defined via the spectral curves
and have little to do with the monopole fields.  The relationship
between the boundary values of the monopole fields and the spectral
curve should expose a more direct link between the monopoles and the
Einstein metric.

\section{Massless monopoles.}  \label{sec:massless}
By the maximum principle on the Higgs field, monopoles with zero mass
are necessarily flat, and hence trivial on hyperbolic space.  Still,
the zero mass limit of hyperbolic monopoles, which by a rescaling
corresponds to the infinite curvature limit of hyperbolic space,
contains intersting features.  This limit was studied in
\cite{AtiYan,AMuMon,JNoZer} for different reasons.

Given a rational function $f:\bbP^1\rightarrow\bbP^1$, one can produce
a curve $C_f\subset\bbP^1\times\bbP^1$ by
\begin{equation} \label{eq:cf}
C_f=\{(w,z)\in\bbP^1\times\bbP^1\ |\ f(z)=\sigma f\sigma(w)\}
\end{equation}
where $\sigma$ is the antipodal map $\sigma(z)=-1/\bar{z}$,
\cite{AtiYan,AMuMon}.  When $f(z)=k'/(z^N-k)$, the curve $C_f$
contains the parameters of a solution of the Yang-Baxter equation
related to the Potts model.  For a degree $N$ rational map $f$, the
curve $C_f$ has the properties
\begin{enumerate}
\item $C_f$ is a curve of bidegree $(N,N)$ on $\bbP^1\times\bbP^1$
\item $C_f$ is a real curve with respect to
$\tau(w,z)=(\hat{z},\hat{w})$ and has no real points
\item $N(D^+-D^-)\sim 0$ where $D^+$ and $D^-$ are the divisors of the
intersection of $C_f$ with $\bbP^1\times\{z_0\}$ and
$\{w_0\}\times\bbP^1$.
\end{enumerate}
The spectral curve $\Sigma$ of a hyperbolic monopole of mass $m$
satisfies conditions 1 and 2 and a modification of condition 3:\\

3*. $(N+2m)(D^+-D^-)\sim 0$ where $D^+$ and $D^-$ are the divisors
of the intersection

\ \ \ \ \ of $\Sigma$ with $\bbP^1\times\{z_0\}$ and
$\{w_0\}\times\bbP^1$.\\

Thus, it is natural to treat the curves $C_f$ as the zero mass limit
of hyperbolic monopoles.  Another way to write (\ref{eq:cf}) is
\[ C_f=\{(w,z)\in\bbP^1\times\bbP^1\ |\ \langle f(w),
f(z)\rangle=0\}\]
where $\langle\cdot,\cdot\rangle$ is the Hermitian metric on $\bbC^2$,
so $C_f$ detects when the subspaces are orthogonal.

By Lemma~\ref{th:real} the holomorphic sphere satisfies $\langle
q(\hat{w}),q(z)\rangle=w^{-k}\psi(w,z)$ or in
other words the zero set of $\psi$, which defines the spectral curve
of the monopole, is given by
\[\Sigma=\{(w,z)\in\bbP^1\times\bbP^1\ |\ \langle q(w),
q(z)\rangle=0\}.\]

Thus we see that the holomorphic sphere $q:\bbP^1\rightarrow\bbP^k$
resembles closely the holomorphic map $f:\bbP^1\rightarrow\bbP^1$.
Moreover, as the mass of the monopole tends to zero, the image of $q$
tends toward being contained in a line in projective space, giving
$f:\bbP^1\rightarrow\bbP^1\hookrightarrow\bbP^k$ in the limit.  We
interpret the image of $q$ to be ``almost'' contained in a line to
mean that the pull-back of the Kahler form under $q$ is close to the
pull-back of the Kahler form when the image lies in a line.  When the
rational map $f$ is given by radial scattering the claim follows from
\cite{NorAsy}, where it is shown that small non-integer mass
hyperbolic monopoles have boundary values perturbed not too far from
the pull-back of the Kahler form on $\bbP^1$ by the rational map.

The particular line in $\bbP^k$ into which the image of $q$ tends is
not significant, since $q$ is only well-defined up to the action of
$U(k+1)$ on its image.  Since every line in $\bbP^k$ is equivalent up
to this action of $U(k+1)$, the image of $q$ can tend to lie inside
any line.

The motivation of \cite{AtiYan,AMuMon} is to find solutions to the
Yang-Baxter equations that use the spectral curve of a monopole and
resemble the curves $C_f$ from the Potts model.  It may be that the
rational map $f(z)=k'/(z^N-k)$ has an analogue $q$ for each mass.  In
general it is hard to find the holomorphic maps $q$ corresponding to
monopoles, however particularly symmetric examples are more accessible
such as those described in Section~\ref{sec:ch2}.

In \cite{JNoZer} it was shown that the rational map associated to a
hyperbolic monopole can be used to construct an explicit solution of a
degenerate form of the Bogomolny equations obtained from the infinite
curvature limit of hyperbolic space.  This explicit solution was
interpreted as an approximate monopole and the curve $C_f$ defined in
(\ref{eq:cf}) naturally arises as a type of spectral curve.  It was
proven that the approximate monopole can flow to a unique genuine
monopole under a heat flow.  This viewpoint may help with the
question: is there a good way to go straight from the holomorphic
sphere $q$ to the monopole field $(A,\Phi)$?  We would hope to
construct from $q$ an approximate monopole and again prove that a
genuine monopole lies nearby.

\section{Axially symmetric monopoles.}  \label{sec:axial}
When the monopole is axially symmetric, the spectral curve $\Sigma$ is
a collection of $k$ $(1,1)$ curves \[\prod_i(w-a_iz)=0.\]
The reality condition on $\Sigma$ implies that for each $i$ there is a
$j$ such that $a_i=\bar{a}_j$.The curve has mass $m$ if
$L^{k+2m}|_{\Sigma}$ is trivial.  With respect to local
trivialisations of $L^{k+2m}$ in neighbourhoods of
$(w,z)=(0,0)$ and $(w,z)=(\infty,\infty)$ a transition
function can be given by $z^{k+2m}w^{-(k+2m)}$.  A
non-vanishing section over $\Sigma$ can be set to be the constant $1$
in a neighbourhood of $(w,z)=(0,0)$.  Along the curve
$w-a_iz=0$, the transition function is
$z_{k+2m}w^{-(k+2m)}=a_i^{-(k+2m)}$ thus $1\mapsto a_i^{-(k+2m)}$
and one condition that this is a global section over $\Sigma$ is that
the sections over each $w-a_iz=0$ agree at
$(w,z)=(\infty,\infty)$.  Thus
\begin{equation}  \label{eq:axial}
     a_i^{k+2m}=a_j^{k+2m}\in\bbR
\end{equation}
for all $i$ and $j$, and since $\bar{a}_i$ is amongst the $a_j$s, the
number $a_i^{k+2m}$ is real.  When $2m$ is not an integer, the
expression $a_i^{k+2m}$ is still uniquely defined.  In general, such
an expression requires the choice of a branch.  In our case, there is
a well-defined branch of $a_i^{k+2m}$ obtained by continuity as the
mass is varied.  When $m=0$, one can still make sense of the spectral
curve of a ``massless'' monopole (see Section~\ref{sec:massless}) and
in this case it is given by the equation $w^k+(-1)^k\alpha^kz^k=0$
for some $\alpha>0$.  Since condition (\ref{eq:axial}) is a discrete
condition on the $a_i$ we can again use continuity in the mass, and
prove that for $m>0$
\begin{equation}  \label{eq:roots}
     a_j=\alpha\exp\left(\frac{2\pi ij}{k+2m}\right),\ j=(1-k)/2,(3-k)/2,
     \ldots,(k-1)/2.
\end{equation}
Notice that, in agreement with Hitchin \cite{HitCon} p.188, that
$a_i^{k+2m}$ is positive when $k$ is odd and negative when $k$ is
even.  This reflects the real structure on the line bundle.

The coefficients of $\prod_j(w-a_jz)$, for $a_j$ defined in
(\ref{eq:roots}) are all non-zero.  This can be seen from the fact
that each symmetric polynomial in the $a_j$'s strictly increases with
the mass since each $a_j(m)$ creeps along the circle towards the
positive real line.  In Lemma~\ref{th:nondeg} it is proven that the
non-degeneracy of the matrix of coefficients of the defining polyomial
of the spectral curve is equivalent to the vanishing of the cohomology
group $H^0(\Sigma,\cO(k,-2))$.

We have given explicit expressions for spectral curves of hyperbolic
monopoles and for boundary values of hyperbolic monopoles.  Here we
give an explicit expression n explicit expression for the field over
$\bbH^3$ of a charge two hyperbolic monopole.  Choose coordinates
$(z,r)$ where $r$ is the hyperbolic distance from the origin to the
point, and $z$ is a holomorphic coordinate on each sphere of constant
distance from the origin.  In order to give a gauge invariant
expresssion it is convenient to use an associated Hermitian metric $H$
defined over $\bbH^3$ which gives the monopole $(A,\Phi)$ in a
non-unitary gauge by \[ A_{\bar{z}}=0,\ A_z=H^{-1}\partial_zH,\
A_r=(1/2)H^{-1}\partial_rH,\ \Phi=(-i/2)H^{-1}\partial_rH.\] The pair
$(A,\Phi)$ satisfies the Bogomolny equation when $H$ satisfies the
nonlinear equation
\begin{equation}  \label{eq:bog}
\partial_r(H^{-1}\partial_rH)+\frac{(1+|z|^2)^2}{\sinh^2(r)}
\partial_{\bar{z}}(H^{-1}\partial_zH)=0.
\end{equation}
For an axially symmetric centred charge 2 hyperbolic monopole, $H$
looks like
\begin{equation}  \label{eq:expr}
H=\frac{1}{D}
\left(\begin{array}{cc}a(r)+2b(r)|z|^2+|z|^4/a(r)&
(1-b(r)^2)^{1/2}(a(r)-1/a(r))\bar{z}^2\\
(1-b(r)^2)^{1/2}(a(r)-1/a(r))z^2&1/a(r)+2b(r)|z|^2+a(r)|z|^4
\end{array}\right)
\end{equation}
for $D=(1+|z|^2)^2-(2-b(r)(a(r)+1/a(r)))|z|^2$.  The functions $a(r)$
and $b(r)$ satisfy a set of non-linear equations derived from putting
(\ref{eq:expr}) into (\ref{eq:bog}).  One explicit solution is given
by
\begin{equation}   \label{eq:sol}
a(r)={\rm sech}(r)=b(r)
\end{equation}
and this gives a mass $1/2$ monopole.  When $a(r)=e^{-2r}$ and
$b(r)\equiv 0$ we get a solution of a degenerate equation much like
(\ref{eq:bog}) which corresponds to a zero mass monopole.

The holomorphic sphere of the monopole arising from (\ref{eq:sol}) is
\[ q(z)=(\frac{1}{\sqrt{2}}(1+z),\frac{i}{\sqrt{2}}(1-z),z^2).\]

It would be desirable to find a one-parameter family of solutions
$a_m(r)$ and $b_m(r)$ depending on the mass $m$, in particular to get
explicit expressions for fractional mass hyperbolic monopoles.

\section{The Vanishing Theorem}
The proof of Theorem~\ref{th:sphere} requires
$H^0(\Sigma,L^{2m}(0,k-2))=0$.  In this section we will prove a more
general vanishing theorem that has further applications.

\begin{thm}[Vanishing Theorem]
\label{thm:vanishing}
If $\Sigma\subset Z$ is the spectral curve of a hyperbolic monopole of
mass $m$ and charge $k$ then
$$
H^0(\Sigma, L^s(k-2, 0)) = 0
$$
for all $1 \leq s \leq 2m+1$.
\end{thm}
\begin{note}
We have that $L^{2m+k}$ restricted to $\Sigma$ is trivial. So
$H^0(\Sigma, L^s(k-2, 0)) = H^0(\Sigma, L^{s-2m-k}(k-2,0)) =
H^0(\Sigma, L^{2m+2-s}(k-2,0))$.  The last equality uses the real
structure and is actually a conjugate linear isomorphism. So it
suffices to prove the theorem for $1 \leq s \leq m+1$.
\end{note}
\begin{note}
The case $s=1$ (and hence also $s = 2m+1$) is elementary since we have 
$H^0(S, \cO(k-1, -1)) = 0$ for any
degree $(k,k)$ curve $\Sigma\subset Z$.
\end{note}

The method of proof of the Theorem is an adaption of
Hitchin  \cite{HitCon} for the Euclidean case. In summary it is as follows.

\noindent(1) Show that the $H^0(\Sigma, L^s(k-2, 0))$ injects into $H^1(Z, 
L^{s-m}\widetilde{E}(-2,0))$.

\noindent(2) Penrose transform to get a solution $u$ of
\begin{equation}
\label{eq:equation}
(\nabla_A^*\nabla_A - 1 + \Phi^*_s\Phi_s) u = 0
\end{equation}
where $\Phi_s = \Phi + i(s-m-1)$
such that $|u(x)| $ decays asympotically as $x \to 0$ like the maximum of
$x^s$ and $x^{2m+2-s}$.

\noindent(3) Transfer $u(x)$ to $\R^4 - \R^2$ where the operator in
\eqref{eq:operator} becomes positive and we can integrate by parts to show that
$u = 0$.

\section{Holomorphic sections of $L$}
Before we begin the proof we need a result about the
space of holomorphic sections of $L^s(k,0)$ over $Z$. Note that
$L^s(k,0)$ extends to the quadric $Q = \bbP^1 \times \bbP^1$
only when $s$ is an integer.  The result we need  says
that if $s$ is not an integer there are no holomorphic sections
of $L^s$ over $Z$
and if $s$ is an integer they are all obtained by restriction
of holomorphic sections of $L^s(k,0)$ over $Q$.  In this
latter case the Kunneth formula tells us that
$$
H^0(Z, L^s(k,0) ) = H^0(Q, L^s(k,0) ) = H^0(\bbP^1, \cO(s+k)) \otimes 
H^0(\bbP^1, \cO(-s)).
$$
We have
\begin{lemma}   \label{th:ext}
 For any non-negative integer $k$
\begin{equation*}
H^0(Z, L^s(k,0))=\begin{cases} 0 & s \neq 0, -1, \dots , -k \\
                              \C^{(s+k)(-s)}    & s = 0, -1, \dots , -k
\end{cases}
\end{equation*}
\end{lemma}
\begin{proof} 
A local section $f$ of $L^s(k,0)=\cO(k+s,-s)$ pulls back to a real
analytic function $\hat{f}$ defined locally on $\C P^3$ satisfying
$\hat f(\alpha u,\beta v)=\alpha^{k+s}\beta^{-s}\hat f(u,v)$ where $u$ and
$v$ lie in $\bbC^2$.  Equivalently $\hat f (\lambda^{1/2}\mu.u,
\lambda^{-1/2}\mu.v) = \lambda^{k/2+s}\mu^k \hat f(u,v)$ and the
factor $\mu^k$ shows that $\hat f$ can be interpreted as a section of
$\cO(k)$ that satisfies
\begin{equation}  \label{eq:pback}
\hat f (\lambda^{1/2}u, \lambda^{-1/2}v) =\lambda^{k/2+s}\hat f(u,v).
\end{equation}
Instead of working on local open neighbourhoods of $\C P^3$ we can
restrict to the set $P$ defined by $P=\{[u,v]\in\C P^3:\langle
u,v\rangle>0\}$ which is an open subset of the real hypersurface
$\{{\rm Im}\langle u,v\rangle=0\}$.  Then (\ref{eq:pback}) describes a
real analytic section of $\cO(k)|_P$, holomorphic on holomorphic
sub-manifolds of $P$, that transforms under $\lambda\in\bbR^+$.  The
set $P$ misses the pull-back of the anti-diagonal and it contains the
pre-images of all real $(1,1)$ curves.  One can describe $P$ as a twistor
space.

Pick any projective line $L \subset P$. Then $\hat f $ continues
analytically to a holomorphic section of $\cO(k)$ in an open (in $\C
P_3$) neighbourhood $W $ of $L$.  We claim that $\hat f$ is the
restriction of a polynomial of degree $k$.  Granted that, it follows
at once that $\hat f = 0$ unless $ s = 0, -1, \dots , -k $ since the
only possible weights for the $\R^+$ action on a polynomial of degree
$k$ are $k/2,k/2-1,k/2-2,...,-k/2$.

It remains to show that $\hat f$ is a restriction of such a 
polynomial.  The identity
$$
\hat f(z_0, z_1, z_2, z_3) = \sum z_{j_1} \dots z_{j_k} 
\partial_{j_1} \dots \partial_{j_k}\hat f
$$
(which follows from Euler's identity proved by repeatedly
differentiating both sides of $f(\lambda z_0,\lambda z_1,\lambda
z_2,\lambda z_3)=\lambda^kf(z_0, z_1, z_2, z_3)$ with respect to
$\lambda$) reduces to the case $k=0$ for $\partial_{j_1} \dots
\partial_{j_k}\hat f \in H^0(W, \cO)$.  But this is constant as $W$
contains lots of intersecting projective lines.
\end{proof}

\section{Proof of the vanishing theorem}
\subsection{The injection}
 From the short exact sequence of sheaves
$$
0 \to \cO(-k, -k) \to \cO \to \cO_\Sigma \to 0
$$
we obtain
$$
0 \to H^0(Z, L^s(k-2,0))   \to H^0(\Sigma, L^s(k-2,0))   \to H^1(Z, 
L^s(-2,-k)) \to \dots .
$$
According to Lemma \ref{th:ext}, $H^0(Z, L^s(k-2,0)) = 0$ 
unless $s=0, -1, -2, \dots 2-k$
so for the range of $s$ we are interested in the connecting homomorphism
$$
H^0(\Sigma, L^s(k-2,0))   \to H^1(Z, L^s(-2,-k))
$$
is injective.

We also have
$$
0 \to L^m(0, -k) \to \widetilde{E} \to L^{-m}(0, k)  \to 0
$$
and hence
$$
0 \to L^{s}(-2, -k) \to L^{s-m}\widetilde{E}(-2,0) \to L^{s-2m}(-2, k)  \to 0
$$
so that
$$
H^0(Z, L^{s-2m}(-2, k)) \to  H^1(Z, L^{s}(-2, -k))  \to H^1(Z, L^{s-m}\widetilde{E}(-2,0)).
$$
Again from Lemma \ref{th:ext}  we have that $H^0(Z, L^{s-2m}(-2, k)) = 
0$ unless $s-2m = 2, 3, \dots, -k$ or
$s = 2m+2, 2m+3, \dots, 2m+k$ which is outside the range of interest. 
Finally we have
$$
H^0(\Sigma, L^s(k-2,0)) \hookrightarrow H^1(Z, L^{s-m}\widetilde{E}(-2,0))
$$
for $1 \leq s \leq 2m+1 $.

By replacing $s $ by $s-2m-k$ at the outset and using the other sequence
$$
0 \to L^{-m}(-k, 0) \to \widetilde{E} \to L^{m}(-k, 0)  \to 0
$$
we obtain an equivalent description of the same class in $H^1(Z, 
L^{s-m}\widetilde{E}(-2,0))$
which factors through the other canonical subbundle.

\subsection{The Penrose Transform for $\bbH^3$}
We describe the Penrose transform of a class
$$
[f] \in H^1(Z, L^a(-2,0))
$$
and give estimates for its growth if it is compactly supported.

Identify hyperbolic three-space $\bbH^3$ with the space of 
positive-definite, two by two,  Hermitian
matices $X$ up to scale or with the space of positive-definite, two 
by two, Hermitian
matrices of unit determinant. We co-ordinatize these matrices by
$$
X = \frac{1}{x_3} \left[ \begin{matrix} 1 & x_1 + ix_2 \\
                                     x_1 -ix_2  & x_3^2 + x_1^2 + x_2^2
                        \end{matrix} \right]
$$
for $x_3 > 0$. This is the upper half-space model of hyperbolic space.
Denote by $M$ the  open set of future pointing timelike vectors in 
$\R^{3,1}$ so that $M/\R^+ = \bbH^3$.

\begin{thm}
\label{th:penrose}
There is a canonical isomorphism
\begin{eqnarray*}
H^1(Z, L^a(-2, 0)) &= \{ u \in C^\infty(M) \colon \square u = 0, Eu = 
(a-2)u \}    \\
                    &= \{ v \in C^\infty(\bbH^3) \colon (\triangle + 
a(a-2))v = 0 \}
\end{eqnarray*}
Moreover if $f$ has compact support then $u$ has a decomposition
for $x \leq 0$ smooth
\begin{equation}
\label{eq:decomposition}
u(x_1,x_2,x_3)  = x_3^{2-a} u_1(x_1,x_2,x_3) + x_3^au_2(x_1,x_2,x_3)
\end{equation}
where $u_1$ and $u_2$ are smooth down to $x_3 = 0$.
\end{thm}
\begin{remark}
  Notice that the expansion fits well with what is known about
the boundary behaviour of eigenfunctions of the Laplacian on $\bbH^3$ 
\cite{MazMel}.
\end{remark}
\begin{proof}
Choose homogeneous coordinates $[\eta]=[\eta^0,\eta^1]$ and
$[\zeta]=[\zeta^0,\zeta^1]$ for $Z$, so in terms of the affine
coordinates, $w=\eta^1/\eta^0$ and $z=\zeta^1/\zeta^0$.  By
considering the framing $\<\eta,\zeta\>/\<\zeta,\zeta\>$ of $L$ over
$Z$, we see that $[f]$ can be represented by
$$
f \in \Omega^{0,1}(Z, \cE(-2,0))
$$
such that
$$
\bar\partial 
\left(\left(\frac{1+\bar{w}z}{1+\bar{z}z}\right)^a f \right) 
= 0.
$$
By pulling back to $P$, for example, we get the function
\begin{equation}
\label{eq:integral}
\psi(x) = \int_{[\zeta] \in \bbP^1}
\left(\frac{\zeta^*X\zeta}{\zeta^*\zeta}\right)^a f(X\zeta, \zeta) \wedge
(\zeta^0d\zeta^1 - \zeta^1 d\zeta^0)
\end{equation}
from the Minkowski version of the Penrose transform. Thus $\psi$ 
satisfies the wave equation
on $M$ and is plainly homogeneous of degree $a-2$ in $X$.  One way to find
the equation satisfied by $\psi$ on restriction to $\bbH^3$ is to
compute $|X|^2 \square ( |X|^{2-a} \psi(X))$, for $\square $ on functions of
degree $0$ is $\triangle_{\bbH^3} = - \tr \Hess$. We have
$$
\partial_i (|X|^{2-a} \psi(X)) = (2-a) |X|^{-a} X_i \psi(X) +
                                          |X|^{2-a} \partial_i \psi(X) .
$$
and
\begin{align*}
\partial^i\partial_i (|X|^{2-a} \psi(X)) &=
(2-a) |X|^{-a} (4\psi(X) - 2 \psi(X) - (2-a) \psi(X)) \\
&= a(2-a)\psi(X) |X|^{-a}.
\end{align*}
So if $v(x) = |X|^{2-a} \psi(X)$ we have $\triangle v + a(a-2) v = 0$.

Now suppose for definiteness that $f$ has support contained in the 
relatively compact subset $V$
defined by the condition
$
  |\eta^0 \zeta^1 - \eta^1 \zeta^0|^2 < R^2 |\eta^0 \bar \zeta^0 + \eta^1\bar \zeta^1|^2
$
that is $| w - z | < R|1 + w \bar{z} |$ for some $R < 0$.
The set $V$ obviously does not intersect the antidiagonal.  Consider 
how it meets a
real curve
$$
C_{(x_1,x_2,x_3)} = \{ \eta^1 \zeta^0 + (x_1+ix_2) \eta^1\zeta^1 -
(x_1-ix_2) \eta^0\zeta^0 - (x_1^2+x_2^2 + x_3^2) \eta^0\zeta^0 = 0 \}
$$
when $x_3$ is small but positive.  Now
$$
C_{(x_1, x_2,0)} = \{ \zeta^0 + (x_1+ix_2) \zeta^1 = 0 \} \cup \{
\eta^1 - (x_1-ix_2)\eta^0 = 0 \}
$$
and
$$
V \cap C_{(0, 0, 0)} = V_1 \cup V_2
$$
where
$$
V_1 = \{\eta^1 = 0 \} \times \{ [\zeta] \colon |\zeta^1/\zeta^0| < R \}
$$
and
$$
V_2 = \{\zeta^0 = 0 \} \times \{ [\eta] \colon |\eta^0/\eta^1| < R \}
$$
with $V_1 \cap V_2 = \emptyset$.
Similarly, for sufficiently small $x_3$,  the intersection of $V$ with
$C_{(x_1,x_2,x_3)}$ is a union of two disjoint sets $V_1$ and $V_2$ which
are approximately of the form
$$
V_1 = \{ \eta^1 = (x_1-ix_2) \eta^0 \} \times U_1
$$
and
$$
V_2 = U_2 \times \{ \zeta^0 + (x_1+ix_2) \zeta^1 = 0 \}
$$
where $U_1$ and $U_2$ are discs.
The decomposition \eqref{eq:decomposition} of $u$ in the statement of
the theorem corresponds to a decomposition of the integral \eqref{eq:integral}
into integrals over $V_1$ and $V_2$.  To see this we check the
growth rates of the two contributions by taking $x_1=x_2=0$. If we integrate
over $V_1$ we obtain
\begin{align*}
&\int_{V_1}\left(\frac{x_3^{-1}|\zeta^1|^2 + x_3|\zeta^1|^2}{|\zeta^0|^2 +
|\zeta^1|^2}\right)^a f(x_3^{-1}\zeta^0, x_3\zeta^1, \zeta^0, \zeta^1)
(\zeta^0 d\zeta^1 - \zeta^1 d\zeta^0)\\ &\quad = x_3^{2-a} \int_{V_1}
\left(\frac{1+x_3^2|z|^2}{1+|z|^2}\right)^a f(1, x_3^2z, 1, z) dz\\
&\quad= x_3^{2-a}u_1(x_3, 0, 0)
\end{align*}
Similarly integrating over $V_2$ gives
\begin{align*}
&\int_{V_2}\left(\frac{x|\eta^0|^2+x^{-1}|\eta^1|^2}{x^2|\eta^0|^2 + 
x^{-2}|\eta^1|^2}\right)^a
f(\eta^0, \eta^1, x\eta^0, x^{-1}\eta^1) (\eta^0 d\eta^1 - \eta^1 d\eta^0)\\
&\quad = x^{a} \int_{V_2} 
\left(\frac{1+x^2/|w|^2}{1+x^4/|w|^2}\right)^a 
f(1/w, 1,
x^2/w, 1) d(1/w)\\
&\quad= x^{a}u_2(x, 0, 0).
\end{align*}
The result follows as the integrals have uniformly compact supports at
$x \to 0$.
\end{proof}
\begin{remark} In particular the `boundary data' $u_1$ and $u_2$ arise
as integrals over the generators of $f$. If we use the real structure
to replace $L^a(-2,0)$ by $L^{-a}(0,-2) = L^a(-2,0)$ the roles of the two
generators are swapped but (of course) the conclusion is the same.
\end{remark}

\subsection{Completion of proof}
Coupling to the bundle $\widetilde{E}$ replaces the differential
equation in Theorem \ref{th:penrose} by the analogous coupled equation
\begin{equation}
\label{eq:operator}
(\nabla_A^*\nabla_A - 1 + \Phi^*_s\Phi_s) u = 0
\end{equation}
where $\Phi_s = \Phi + i(s-m-1)$.  The methods of Hitchin \cite{HitCon}
can be used to obtain additional decay like $x_3^m$. So we have
$$
|u(x_3)| \simeq \max( x_3^s , x_3^{2+2m-s}).
$$

It is not clear that the operator in \eqref{eq:operator} is positive so we now
transfer to $\R^4 - \R^2$.  Using $\hat\nabla $ for the Euclidean 
operators we have
$$
x_3^3 \hat\nabla^*\hat\nabla(x_3^{-1}u) = (\nabla^*\nabla-1)u
$$
so that if $v = x_3^{-1}u$ then we have
$$
\hat\nabla^*\hat\nabla v + x_3^{-2} \Phi^*_s\Phi_s v = 0.
$$
In now $1 < s < 2m+1$ then $v = O(x_3^\epsilon)$ for $\epsilon > 0$ and 
this is enough to
integrate by parts in $\R^4 - \R^2$ and hence $v = 0$.

\end{document}